# EDGE PERCOLATION ON A RANDOM REGULAR GRAPH OF LOW DEGREE[1]

BY BORIS PITTEL

*Ohio State University*

*Dedicated to Jimmy and Katie Lyons*

Consider a uniformly random regular graph of a fixed degree $d \geq 3$, with $n$ vertices. Suppose that each edge is *open* (*closed*), with probability $p(q = 1 - p)$, respectively. In 2004 Alon, Benjamini and Stacey proved that $p^* = (d-1)^{-1}$ is the threshold probability for emergence of a giant component in the subgraph formed by the open edges. In this paper we show that the transition window around $p^*$ has width roughly of order $n^{-1/3}$. More precisely, suppose that $p = p(n)$ is such that $\omega := n^{1/3}|p - p^*| \to \infty$. If $p < p^*$, then with high probability (whp) the largest component has $O((p - p^*)^{-2} \log n)$ vertices. If $p > p^*$, and $\log \omega \gg \log \log n$, then whp the largest component has about $n(1 - (p\pi + q)^d) \asymp n(p - p^*)$ vertices, and the second largest component is of size $(p - p^*)^{-2}(\log n)^{1+o(1)}$, at most, where $\pi = (p\pi + q)^{d-1}, \pi \in (0, 1)$. If $\omega$ is merely polylogarithmic in $n$, then whp the largest component contains $n^{2/3+o(1)}$ vertices.

**1. Introduction and results.** Let $d \geq 3$ and $n > d$ be given. Assuming that $nd$ is even, introduce $\mathcal{G}_n$, the sample space of all $d$-regular graphs on the vertex set $[n]$. It is known (Bender and Canfield [5]) that

$$(1.1) \qquad |\mathcal{G}_n| = \frac{(nd - 1)!!}{(d!)^n} \exp\left(-\frac{d^2 - 1}{4} + O(n^{-1})\right).$$

Here is a probabilistic interpretation of (1.1), due to Bollobás [7], who actually considered a general degree sequence. Introduce $n$ sets $S_1, \ldots, S_n$, $S_i$ representing vertex $i$, with $|S_i| \equiv d$. Then $(nd - 1)!!$ is the total number of complete matchings (pairings) on $S = \bigcup_i S_i$, and the exponential factor is

Received July 2006; revised July 2006.
[1]Supported in part by NSF Grant DMS-04-06024.

the probability that the matching chosen uniformly at random is graph-induced, that is, there are no pairs $(s,t)$ with $s$ and $t$ from the same set $S_i$, $i \in [n]$, and no pairs $(s,t)$, $(s',t')$ with $s,s' \in S_i$, $t,t' \in S_j$, $i \neq j$. Equation (1.1) follows then from the observation that every $d$-regular simple graph induces $(d!)^n$ distinct pairings on $S$.

Consider $G_n$ the uniformly random $d$-regular graph on $n$ vertices. Thus $G_n$ assumes every one of its $|\mathcal{G}_n|$ values with the same probability $|\mathcal{G}_n|^{-1}$.

Let $p \in (0,1)$ be fixed. We define the random graph $G_{np}$ as the random subgraph of $G_n$ obtained by "opening" ("closing") each edge of $G_n$ with probability $p$ ($q = 1-p$), independently of all other edges. Thus the edge set of $G_{np}$ is the set of all open edges of $G_n$. Several years ago Itai Benjamini [6] posed a problem (i) to show that $p^* = (d-1)^{-1}$ is a threshold value of $p$ for emergence of a giant component in $G_{np}$, and (ii) to determine the width of the transition window around $p^*$. In 2004 the part (i) was solved by Alon, Benjamini and Stacey [2].

Here is a quick-and-dirty argument for why $p^*$ had better be the threshold. A subgraph of $G_n$ induced by the vertices within a relatively small distance from any given vertex $v$ is a tree in which all nonleaves have degree $d$. Orienting the edges away from $v$, we get a directed tree in which the root $v$ has outdegree $d$, and the remaining nonleaves each have outdegree $d-1$. And it is known (Durrett [13]; see also Grimmett [17]) that, for the edge (bond) percolation on an infinite directed tree of outdegree $d-1$, the critical probability is $p^* = (d-1)^{-1}$, and for $p > p^*$ the probability that the root is in an infinite cluster is $1 - \pi$,

$$(1.2) \qquad \pi = (\pi p + q)^{d-1}, \qquad \pi \in (0,1).$$

Here $\pi = \pi(p)$ is the probability of eventual extinction for the branching process with immediate family size having binomial distribution $\text{Bin}(d-1, p)$. Thus we should anticipate that each of the $d$ neighbors of vertex $v$ will have only few "descendants" with probability $\pi$, independently of other neighbors, and so $v$ itself will have only few descendants with probability

$$\sum_{j=0}^{d} \binom{d}{j} p^j q^{d-j} \pi^j = (p\pi + q)^d.$$

Guessing that with high probability (whp), that is, with probability $1 - o(1)$, each vertex $v$ is either in a small component or in a unique "giant" component, we can even anticipate then that the fraction of vertices in the giant component is asymptotic to

$$(1.3) \qquad \alpha(p) = 1 - (p\pi + q)^d.$$

(This was not proved in [2].)

What about the window width? Since the pioneering work of Erdős and Rényi [16] in the early 1960s, it has been known that for $d = n - 1$, that is, $G_n = K_n$, a unique giant component appears when the average vertex degree $\overline{d}$ exceeds 1. The issue of the transition window around 1 remained wide open until Bollobás [9] (see also [10]) was able to show that its width is of order $O(n^{-1/3}(\log n)^{2/3})$. Later Łuczak [19] extended the techniques of [9, 10] and showed, among other sharp results, that the logarithmic factor could be removed. (See also [18, 21, 25] concerning the critical phase $\overline{d} - 1 = O(n^{-1/3})$.) In light of these results Itai Benjamini wondered whether the bound $n^{-1/3}$ would hold for the graph $G_n$ of a small fixed vertex degree $d$ as well. Recently Borgs et al. [11, 12] have established this width bound for a class of deterministic regular graphs with vertex-transitivity property, meeting a certain "triangle condition." This class contains $K_n$ and some "high-dimensional" tori. Notice that whp the random $d$-regular graph $G_n$ is not in this class, since almost all $d$-regular graphs are asymmetric (Bollobás [8]), whence intransitive.

Our goal is to solve part (ii) of Benjamini's problem, confirming his $n^{-1/3}$ conjecture. Here are our main results.

Given $p$, let $L_n^{(1)}, L_n^{(2)}$ denote the size of the largest component and the size of the second largest component in $G_{np}$.

THEOREM 1 (Subcritical case). *Suppose that $p = p(n)$ is such that*

$$\lim_{n \to \infty} n^{1/3}(p^* - p) = \infty. \tag{1.5}$$

*Then $\frac{L_n^{(1)}}{(p^*-p)^{-2} \log n}$ is bounded in probability, or $L_n^{(1)} = O_P((p^* - p)^{-2} \log n)$ in short.*

THEOREM 2 (Nearcritical case). *Suppose that, for some $a > 0$,*

$$n^{1/3}|p - p^*| \leq (\log n)^a. \tag{1.6}$$

*Then, for any $b > \max\{a, 1/3\}$,*

$$\lim_{n \to \infty} P(n^{2/3}(\log n)^{-2a} \leq L_n^{(1)} \leq n^{2/3}(\log n)^b) = 1. \tag{1.7}$$

THEOREM 3 (Supercritical case). *Suppose that*

$$n^{1/3}(p - p^*) \geq (\log n)^{a(n)}, \tag{1.8}$$

*where $a(n) \to \infty$ however slowly. Then*

$$\begin{aligned} L_n^{(1)} &= (1 + o_P(1))n(1 - (p\pi + q)^d) \asymp_P n(p - p^*), \\ L_n^{(2)} &= (p - p^*)^{-2}(\log n)^{1+o_P(1)}; \end{aligned} \tag{1.9}$$

here $o_P(1)$ *denotes a random variable converging to* 0 *in probability, and the second equality for* $L_n^{(1)}$ *means that, in probability,* $L_n^{(1)}$ *is of order* $n(p-p^*)$ *exactly.*

Asaf Nachmias and Yuval Peres (personal communication) informed me that they have obtained similar (but somewhat more precise) results.

These estimates are very similar to, but not as sharp as, those in Bollobás [9, 10] and Łuczak [19] for $G_n = K_n$. In both papers the proofs relied heavily on Bollobás' striking bound

$$(1.10) \qquad C(\nu, \nu + \ell) \leq (c/\ell)^{\ell/2} \nu^{\nu + (3\ell - 1)/2},$$

$C(\nu, \nu + \ell)$ being the total number of connected graphs with $\nu$ vertices and $\nu + \ell$ edges. (See Łuczak [20] for a short proof of (1.10) with the best constant $c$.) While there is a formula for the number of trees with a given degree sequence (Moon [24]) no analogue of the general bound (1.10) for the number of connected graphs with a given degree sequence seems to be known. Fortunately it is possible to obtain sufficiently sharp estimates for the parameters of a randomized algorithm that determines a graph component containing a given vertex. These estimates, together with an asymptotic formula for the likely count of certain *tree* components in the case $p = p^*$, yield the probabilistic bounds for $L_n^{(1)}$ and $L_n^{(2)}$. More specifically, we derive the equations for the conditional (one-step) expectations of this process, which suggest a system of deterministic differential equations whose solution should be a good approximation for the scaled random trajectory. That this is indeed the case is proved by using exponential (super)martingales constructed from the integrals of these differential equations. More general, martingale-based, techniques for random graphs developed by Wormald [27] do not provide estimates strong enough for our purposes. Our method is reminiscent of a technique used in [26] ($k$-core problem), [3] (Karp–Sipser matching algorithm), and [1] (random graphs with immigrating vertices). But most closely related is our joint paper with Jozsi Balogh [4] in which we used the differential equations and the exponential supermartingales for analysis of a "bootstrap" site percolation on the random $d$-regular graph. We must also mention that Molloy and Reed [22, 23] had used the differential equations in their study of existence of a giant component in a uniformly random graph with a given degree sequence. A more or less direct application of their remarkable result to $G_{np}$ seems fraught with subtle difficulties, even for $p$ bounded away from $p^*$. Certainly our bounds for the transition window width could not be obtained this way. However, an algorithmic approach to the giant component phenomenon ushered in [22, 23] is a key tool in both [4] and the present paper.

We believe that our techniques can be extended to a more general (uniformly) random graph $G_{n,\mathbf{d}}$ with a given degree sequence $\mathbf{d} = (d_1, \ldots, d_n)$, such that each $d_i \in [3, d]$, $d < \infty$. Intuitively, in this case the counterpart of (1.2) is

$$\pi = \sum_{j=3}^{d} \lambda_j (\pi p + q)^{j-1}, \qquad \lambda_j := \frac{\sum_{i: d_i = j} d_i}{\sum_{i' \in [n]} d_{i'}}$$

(can the reader see why?), and so

$$p^* = \frac{1}{\sum_j j \lambda_j - 1}.$$

And (1.3) becomes

$$\alpha(p) = 1 - \sum_j \lambda_j (\pi p + q)^j.$$

We think that again the transition window around $p^*$ has width of order $n^{-1/3}$.

In conclusion we should mention a new book by Durrett [15] in which the random graph processes are used systematically either as proof tools or as probabilistic models of evolving networks. The percolation on $G_{n,\mathbf{d}}$ we have just described is inspired by the discussion of the largest cluster problem for a random graph, with a given degree distribution, in Chapter 1 of this book.

The rest of the paper is organized as follows. In Section 2 ($p$ close to $p^*$) we evaluate sharply two first moments of the number of tree components of sizes dependent on $|p - p^*|$, and deduce the lower bound for $L_n^{(1)}$, that is, the left half of (1.7) in Theorem 2. In Section 3 we introduce the percolation process on the random *pairing* and show that it is enough to analyze this eminently more tractable process. Motivated by the formulas for the conditional expected state changes of the attendant Markov chain, we introduce and solve the deterministic system of differential equations, and prove—via a family of the exponential supermartingales—that whp the random process stays sufficiently close to the deterministic trajectory. In Sections 5, 6 and 7 we use the probability tail estimates to complete the proof of the Theorem 2, and to prove Theorems 1 and 3.

**2. Proof of Theorem 2 (Nearcritical case). Lower bound.** We begin at the middle since the proof is purely enumerational, a natural extension of an argument for $K_n$ by Bollobás [9, 10].

First of all, from the discussion of (1.1) in the Introduction we know that, for the uniformly random pairing $\mathcal{P}_n$ on the set $S = \bigcup_{i=1}^{n} S_i$, and the event $A_n = \{\mathcal{P}_n \text{ is simple graph-induced}\}$,

(2.1) $\quad \mathrm{P}(A_n) = \exp\left(-\frac{d^2 - 1}{4} + O(n^{-1})\right) \to \exp\left(-\frac{d^2 - 1}{4}\right) > 0.$

So we will study percolation on the random pairing $\mathcal{P}_n$ instead of percolation on the graph $G_n$. Here we open each *pair* $(s,t) \in \mathcal{P}_n$ with probability $p$, independently of other pairs. And the problem becomes: find the probabilistic bounds for the sizes of the largest component and the second largest component of the (multi)graph $MG_{np}$ on $[n]$ determined by the subpairing $\mathcal{P}_{np}$ consisting of the open pairs. [Each $(s,t) \in \mathcal{P}_n$, $s \in S_i$, $t \in S_j$, gives rise to an edge $(i,j)$ if $i \neq j$, a loop at $i$ if $i = j$.] According to (2.1), for every set $\mathcal{H}$ of graphs on $[n]$,

$$(2.2) \qquad \mathrm{P}(G_{np} \in \mathcal{H}) = \mathrm{P}(MG_{np} \in \mathcal{H} | A_n) = O(\mathrm{P}(MG_{np} \in \mathcal{H})).$$

*Thus any rare event for the random multigraph $MG_{np}$ is a rare event for the random graph $G_{np}$.*

Let

$$n^{1/3}|p - p^*| \leq (\log n)^a.$$

We need to show that whp

$$L_n^{(1)} \geq n^{2/3}(\log n)^{-2a}.$$

It is enough to prove that whp the size of the largest *tree* component exceeds $n^{2/3}(\log n)^{-2a}$.

For $k \geq 1$, let $E_k$ denote the expected number of tree components of $MG_{np}$ with $k$ vertices, and for $k^1, k^2 \geq 1$ let $E_{k^1,k^2}$ denote the expected number of ordered pairs of tree components, with $k^1$ and $k^2$ vertices, respectively.

LEMMA 1. (i) *If $k = o(n^{2/3})$, then*

$$(2.3) \qquad E_k = ndp^{k-1}q^{kd-2(k-1)} \frac{(k(d-1))!}{k!(kd-2(k-1))!} e^{\Psi_{nk}},$$

$$(2.4) \qquad \Psi_{nk} := \alpha \frac{k^2}{n} + o(1),$$

$$(2.5) \qquad \alpha := \frac{(d-1)(d-2)}{2dq}(p - p^*), \qquad p^* = \frac{1}{d-1}.$$

(ii) *If $k^1, k^2 = o(n^{2/3})$, then*

$$(2.6) \qquad \frac{E_{k^1,k^2}}{E_{k^1} E_{k^2}} = 1 + \alpha \left( \frac{k^2}{n} - \sum_{t=1}^{2} \frac{(k^t)^2}{n} \right) + o(1), \qquad k = k^1 + k^2.$$

NOTES. (1) Observe that in (2.5) $\alpha = 0$ for $p = p^*$. (2) For $k = 1$, we obtain

$$E_1 \sim nq^d, \qquad E_{1,1} \sim (nq^d)^2.$$

So the number of isolated vertices, with no loops, is whp asymptotic to $nq^d$, which should be expected, of course.

PROOF OF LEMMA 1. (i) By symmetry,

$$E_k = \binom{n}{k} P_k, \tag{2.7}$$

where $P_k$ is the probability that $[k] = \{1, \ldots, k\}$ is a vertex set of a tree component. Let us compute this probability. Let $\boldsymbol{\delta} = (\delta_1, \ldots, \delta_k)$ denote a generic degree sequence of a tree on $[k]$. It is well known that $\boldsymbol{\delta} > \mathbf{0}$ is a degree sequence of such a tree iff

$$\sum_{i \in [k]} \delta_i = 2(k-1),$$

and that the number of trees with a given feasible $\boldsymbol{\delta}$ is [24]

$$\frac{(k-2)!}{\prod_{i \in [k]}(\delta_i - 1)!}.$$

Besides, each $\delta_i$ is at most $d$, since we are dealing with a subgraph of $G_n$. For a corresponding subpairing of $\mathcal{P}_n$ we need to pick $\delta_i$ points from the set $S_i$, $i \in [k]$. Once it is done, the number of tree-induced subpairings on the chosen points is

$$\frac{(k-2)!}{\prod_{i \in [k]}(\delta_i - 1)!} \cdot \prod_{i \in [k]} \delta_i! = (k-2)! \prod_{i \in [k]} \delta_i.$$

Let $\Delta_1, \ldots, \Delta_k$ be the generic numbers of additional points from the sets $S_1, \ldots, S_k$ which are paired among themselves rather than with points in the sets $S_{k+1}, \ldots, S_n$. It is necessary, of course, that $\delta_i + \Delta_i \leq d$, $i \in [k]$, and that $\Delta := \sum_{i \in [k]} \Delta_i$ is even. The number of ways to pair these $\Delta$ points is $(\Delta - 1)!!$. The remaining $\partial_i = d - \delta_i - \Delta_i$ points in each set $S_i$, $\partial = \sum_{i \in [k]} \partial_i$ in total, are to be paired with some $\partial$ points from $\bigcup_{i \notin [k]} S_i$, while the rest of points from $\bigcup_{i \notin [k]} S_i$, all $(n-k)d - \partial$ of them, are to be paired among themselves. Combined, this can be done in

$$\binom{nd - kd}{\partial} \partial!((n-k)d - \partial - 1)!!$$

ways. In addition, the number of ways to split each $S_i$ into three ordered subsets of cardinalities $\delta_i, \Delta_i, \partial_i$ is

$$\prod_{i \in [k]} \binom{d}{\delta_i, \Delta_i, \partial_i}.$$

The probability that a tree subgraph, with the degree parameters $\boldsymbol{\delta}, \boldsymbol{\Delta}, \boldsymbol{\partial}$, is a component of the random sub(multi)graph $MG_{np}$ is

$$p^{k-1} q^{\Delta/2 + \partial}.$$

Let $P(\boldsymbol{\delta}, \boldsymbol{\Delta}, \boldsymbol{\partial})$ denote the probability that $G_{np}$ has a tree component with vertex set $[k]$, and the degree parameters $\boldsymbol{\delta}, \boldsymbol{\Delta}, \boldsymbol{\partial}$. Collecting all the pieces, we obtain

$$
\begin{aligned}
P(\boldsymbol{\delta}, \boldsymbol{\Delta}, \boldsymbol{\partial}) \\
(2.8) \quad &= \frac{1}{(nd-1)!!}(k-2)!(\Delta-1)!!\binom{(n-k)d}{\partial}\partial!((n-k)d-\partial-1)!! \\
&\quad \times \prod_{i \in [k]} \delta_i \, d\, \delta_i, \Delta_i, \partial_i \cdot p^{k-1} q^{\Delta/2+\partial}.
\end{aligned}
$$

Notice that $\Delta$ uniquely determines $\partial$,

$$\partial = kd - (k-1) - \Delta.$$

Let $P(\Delta)$ stand for the probability of a tree component with a given $\Delta$; so $P(\Delta)$ is the sum of $P(\boldsymbol{\delta}, \boldsymbol{\Delta}, \boldsymbol{\partial})$'s over $\boldsymbol{\delta}$'s, $\boldsymbol{\Delta}$'s and $\boldsymbol{\partial}$'s satisfying the conditions

$$
(2.9) \quad \begin{aligned}
\Delta_i + \delta_i + \partial_i &= d, \quad i \in [k], \\
\sum_{i \in [k]} \delta_i &= 2(k-1), \\
\sum_{i \in [k]} \Delta_i &= \Delta.
\end{aligned}
$$

According to (2.8), given $\Delta$, we need to evaluate

$$
\begin{aligned}
\sum_{\substack{\boldsymbol{\Delta}, \boldsymbol{\delta}, \boldsymbol{\partial} \\ \text{meet (2.9)}}} \prod_{i \in [k]} & \frac{1}{(\delta_i-1)!\Delta_i!\partial_i!} \\
= [x^\Delta y^\partial] & \sum_{\substack{\mathbf{a}, \boldsymbol{\delta}, \boldsymbol{\partial} \geq \mathbf{0} \\ a_i + \delta_i + \partial_i = d-1, i \in [k]}} \prod_{i \in [k]} \frac{1}{a_i! \Delta_i! \partial_i!} \frac{x^{\Delta_i} y^{\partial_i}}{} \\
= & \frac{1}{[(d-1)!]^k}[x^\Delta y^\partial] \prod_{i \in [k]} \sum_{\substack{a_i, \Delta_i, \partial_i \geq 0 \\ a_i + \Delta_i + \partial_i = d-1, i \in [k]}} \binom{d-1}{a_i, \Delta_i, \partial_i} x^{\Delta_i} y^{\partial_i} \\
= & \frac{1}{[(d-1)!]^k}[x^\Delta y^\partial]((1+x+y)^{(d-1)})^k \\
= & \frac{1}{[(d-1)!]^k} \cdot \frac{(k(d-1))!}{(k-2)!\Delta!\partial!}.
\end{aligned}
$$

Therefore, summing $P(\boldsymbol{\delta},\boldsymbol{\Delta},\boldsymbol{\partial})$'s over all $\boldsymbol{\delta}$'s, $\boldsymbol{\Delta}$'s, $\boldsymbol{\partial}$'s for a given $\Delta$, we obtain

$$P(\Delta) = \frac{d^k(k(d-1))!}{(nd-1)!!} \frac{(\Delta-1)!!}{\Delta!} \binom{(n-k)d}{\partial} ((n-k)d - \partial - 1)!!$$

$$(2.10) \qquad \times p^{k-1} q^{\Delta/2 + \partial},$$

$$\partial = kd - 2(k-1) - \Delta.$$

Using $(2s-1)!!/(2s)! = (2^s s!)^{-1}$, substituting $\Delta = 2D$ and excluding $\partial$, we rewrite (2.10) as follows:

$$P(2D) = \frac{d^k 2^{kd-k+1} (k(d-1))! ((n-k)d)!}{2^{nd/2}(nd-1)!!} \cdot p^{k-1} q^{kd-2k+2}$$

$$(2.11)$$

$$\times \frac{(4q)^{-D}}{D!(kd-2k+2-2D)!(nd/2 - kd + k - 1 + D)!}.$$

What is left is to find a sharp estimate for the sum of $P(2D)$ over $D \leq kd - 2(k-1)$. Notice that the $D$-dependent fraction in (2.11), call it $F(D)$, is bounded above by

$$F^+(D) = \frac{1}{(kd-2k+2)!(nd/2 - kd + k - 1)!} \cdot \frac{1}{D!} \left( \frac{(kd-2(k-1))^2}{4q(nd/2 - kd + k - 1)} \right)^D.$$

Introduce a Poisson($\lambda_k$) random variable $Z$, where

$$(2.12) \quad \lambda_k = \frac{(kd-2(k-1))^2}{4q(nd/2 - kd + k - 1)} = \frac{k^2(d-2)^2}{2qnd} + O(k/n + k^3/n^2).$$

Observe that $\lambda_k$ is of order $k^2/n$ exactly. Recall that $k^3/n^2 \to 0$. Using Chebyshev's inequality,

$$P\left( Z \geq \lambda_k + \lambda_k^{1/2} \left(\frac{n^2}{k^3}\right)^{\varepsilon} \right) \leq \left(\frac{k^3}{n^2}\right)^{2\varepsilon} \to 0, \qquad n \to \infty,$$

for every fixed $\varepsilon > 0$. Therefore, denoting $\lambda_k + \lambda_k^{1/2}(n^2/k^3)^{\varepsilon}$ by $D_k$,

$$(2.13) \quad \sum_{D \geq D_k} F^+(D) \leq \frac{e^{\lambda_k}}{(kd-2k+2)!(nd/2 - kd + k - 1)!} \cdot \left(\frac{k^3}{n^2}\right)^{2\varepsilon}.$$

For $D \leq D_k$, we have

$$[(kd-2k+2-2D)!]^{-1} = [(kd-2k+2k)!]^{-1} \prod_{j=1}^{2D} (kd-2k+2-j)$$

$$= [(kd-2k+2)!]^{-1}(kd-2k+2)^{2D} e^{O(D_k^2/k)},$$

where
$$\frac{D_k^2}{k} \leq 2\frac{\lambda_k^2}{k} + 2\frac{\lambda_k}{k}\left(\frac{n^2}{k^3}\right)^{2\varepsilon} = O\left(\frac{k^3}{n^2} + n^{-1/3}\left(\frac{k}{n^{2/3}}\right)^{1-6\varepsilon}\right) \to 0,$$

provided that $\varepsilon = 1/6$, say. So
$$[(kd - 2k + 2 - 2D)!]^{-1} = (1 + o(1))[(kd - 2k + 2)!]^{-1}(kd - 2k + 2)^{2D},$$

and likewise
$$[(nd/2 - kd + k - 1 + D)!]^{-1} = (1 + o(1))[(nd/2 - kd + k - 1)!]^{-1}.$$

Consequently
$$\sum_{D \leq D_k} F(D) = (1 + o(1)) \sum_{D \leq D_k} F^+(D)$$

$$(2.14) \qquad = (1 + o(1))\frac{e^{\lambda_k}}{(kd - 2k + 2)!(nd/2 - kd + k - 1)!}\mathrm{P}(Z < D_k)$$

$$= (1 + o(1))\frac{e^{\lambda_k}}{(kd - 2k + 2)!(nd/2 - kd + k - 1)!}.$$

Combining (2.11), (2.13) and (2.14), we get
$$P_k = \sum_{D \leq (kd-2k+2)/2} P(2D)$$

$$(2.15)$$
$$= (1 + o(1))\frac{d^k 2^{kd-k+1}(k(d-1))!((n-k)d)!p^{k-1}q^{kd-2k+2}}{2^{nd/2}(nd-1)!!(kd - 2k + 2)!(nd/2 - kd + k - 1)!} \cdot e^{\lambda_k}.$$

Call the last fraction $\mathcal{F}_k$. $\mathcal{F}_{nk}$, the $n$-dependent portion of $\mathcal{F}_k$, can be rewritten this way:

$$\mathcal{F}_{nk} = \frac{((n-k)d)!}{2^{nd/2}(nd-1)!!((n/2-k)d+(k-1))!}$$

$$= \frac{((n-k)d)!}{((n-k)d+k-1)!} \cdot \frac{\binom{(n-k)d+k-1}{\frac{nd}{2}}}{\binom{nd}{\frac{nd}{2}}}.$$

Evaluating each of the two factors on the right in a standard way, we obtain
$$\mathcal{F}_{nk} = (nd)^{-k+1} 2^{-kd+k-1} e^{\Phi_{nk}},$$

where
$$(2.16) \qquad \Phi_{nk} := \frac{k^2}{n}\left(1 - \frac{1}{2d} - \frac{(d-1)^2}{2d}\right) + O(k/n + k^3/n^2).$$

Therefore the formula for $\mathcal{F}_k$ simplifies to

$$(2.17) \qquad \mathcal{F}_k = n^{-(k-1)} d p^{k-1} q^{kd-2(k-1)} \frac{((k(d-1))!}{(kd-2(k-1))!} e^{\Phi_{nk}}.$$

Using (2.7), (2.15), together with (2.12), (2.16), (2.17), and

$$\binom{n}{k} = \frac{n^k}{k!} \exp\left(-\frac{k^2}{2n} + O(k/n + k^3/n^2)\right),$$

we obtain then

$$(2.18) \qquad E_k = \binom{n}{k} P_k = n d p^{k-1} q^{kd-2(k-1)} \frac{(k(d-1))!}{k!(kd-2(k-1))!} e^{\Psi_{nk}},$$

where

$$\Psi_{nk} = \lambda_k + \Phi_{nk} - \frac{k^2}{2n} + o(1)$$
$$= \alpha \frac{k^2}{n} + o(1),$$
$$\alpha := \frac{(d-1)(d-2)}{2dq}(p - p^*), \qquad p^* = \frac{1}{d-1}.$$

The proof of part (i) is complete.

Part (ii). Since the proof parallels the above argument, we will skip some details. First of all,

$$(2.19) \qquad E_{k^1,k^2} = \binom{n}{k^1, k^2, n-k} P_{k^1,k^2}, \qquad k := k^1 + k^2,$$

where $P_{k^1,k^2}$ is the probability that the sets $[k^1] = \{1,\ldots,k^1\}$, $[k^2] := \{k^1 + 1,\ldots,k\}$ are the vertex sets of the tree components $T^1$ and $T^2$ in $MG_{np}$. Let $(\boldsymbol{\delta}^1, \boldsymbol{\Delta}^1, \boldsymbol{\partial}^1)$ and $(\boldsymbol{\delta}^2, \boldsymbol{\Delta}^2, \boldsymbol{\partial}^2)$ denote the degree parameters of $T^1$ and $T^2$, as in Part (i). Then $P(\boldsymbol{\delta}, \boldsymbol{\Delta}, \boldsymbol{\partial})$, the probability that $MG_{np}$ has the tree components $T^1$ and $T^2$ with those parameters, is given by

$$P(\boldsymbol{\delta}, \boldsymbol{\Delta}, \boldsymbol{\partial})$$
$$= \frac{1}{(nd-1)!!} \prod_{s=1}^{2} (k^s - 2)!(\Delta^s - 1)!!$$
$$\times \prod_{i \in [k^s]} \delta_i^s \binom{d}{\delta_i^s, \Delta_i^s, \partial_i^s} p^{k^s-1} q^{\Delta^s/2}$$
$$(2.20) \qquad \times \sum_a q^{\partial - a} \binom{\partial^1}{a} \binom{\partial^2}{a} a!$$

$$\times \binom{(n-k)d}{\partial^1 - a, \partial^2 - a, (n-k)d - \partial + 2a} \prod_{t=1}^{2}(\partial^t - a)!$$

$$\times ((n-k)d - \partial + 2a - 1)!!,$$

where $\partial = \partial^1 + \partial^2$. The first line expression contains a double dose of some of the corresponding factors in (2.8). In particular,

$$\prod_{s=1}^{2} p^{k^s-1} q^{\Delta^s/2}$$

is the probability that all the edges of $T^1$ and $T^2$ are open, and all other edges induced by the vertex sets $[k^1]$ and $[k^2]$ are closed. As for the second and third lines sum, the $a$th summand is the number of ways (1) to choose $a$ points among $\partial^1$ points from $\bigcup_{i \in [k^1]} S_i$ and $a$ points among $\partial^2$ points from $\bigcup_{i \in [k^2]} S_i$ and match them, then (2) to choose among the remaining $(n-k)d$ points $\partial^t - a$ partners for the remaining $\partial^t - a$ points from $\bigcup_{i \in [k^t]} S_i$, $t=1,2$, and finally (3) to match the rest of the $(n-k)d - \partial + 2a$ points among themselves, multiplied by the probability $q^{\partial - a}$ that all the pairs $(u,v)$, such that $u \in \bigcup_{i \in [k^t]} S_i$, $v \notin \bigcup_{i \in [k^t]} S_i$, $t=1,2$, are closed.

Given numbers $D^t$, let $P(2D^1, 2D^2)$ denote the probability that $MG_{np}$ has these tree components $T^1$ and $T^2$ with the parameters $\Delta^t = 2D^t$, $t=1,2$. Summing $P(\boldsymbol{\delta}, \boldsymbol{\Delta}, \boldsymbol{\partial})$ over feasible $\boldsymbol{\delta}$'s, $\boldsymbol{\Delta}$'s and $\boldsymbol{\partial}$'s, we obtain

$$P(2D^1, 2D^2)$$

$$= \frac{d^k((n-k)d)! 2^{kd-k+2} \prod_{s=1}^{2}(k^s(d-1))! p^{k-2} q^{kd-2k+4}}{2^{nd/2}(nd-1)!!}$$

(2.21) $$\times \mathcal{S}(D^1, D^2),$$

$$\mathcal{S}(D^1, D^2)$$

$$= \sum_{a} \frac{2^{-2D-a} q^{-D-a}}{a! \prod_{t=1}^{2} D^t!(k^t(d-2) + 2 - a - 2D^t)!(nd/2 - kd + k - 2 + a + D)!},$$

where $D = D^1 + D^2$. Acting as in part (i), we upperbound the $a$th summand in (2.21) by

$$\frac{1}{a!} \left( \frac{\prod_{t=1}^{2}(k^t d - 2k^t + 2)}{q(nd - 2kd + 2k - 4)} \right)^a \cdot \prod_{\tau=1}^{2} \frac{1}{D^\tau!} \left( \frac{(k^\tau d - 2k^\tau + 2)^2}{2q(nd - 2kd + 2k - 4)} \right)^{D_\tau}$$

$$\times \left( \prod_{t=1}^{2}(k^t(d-2) + 2)! \right)^{-1} ((nd/2 - kd + k - 2)!)^{-1}.$$

Then the sum of the first line product over $a$, $D^1$, $D^2$ is at most

$$\exp\left[\frac{\prod_{t=1}^2(k^t(d-2)+2)}{q(nd-2kd+2k-4)} + \sum_{t=1}^2 \frac{(k^t(d-2)+2)^2}{2q(nd-2kd+2k-4)}\right]$$

$$= \exp\left[\frac{\prod_{t=1}^2 k^t(d-2)}{qnd} + \sum_{t=1}^2 \frac{(k^t(d-2))^2}{2qnd} + O(k/n + k^3/n^2)\right]$$

$$= \exp\left[\frac{(k(d-2))^2}{2qnd} + O(k/n + k^3/n^2)\right], \qquad k = k^1 + k^2.$$

And, in fact, that sum is asymptotic to the last expression. This can be proven by using three independent Poissons, with parameters

$$\lambda^t = \frac{(k^t(d-2)+2)^2}{2q(nd-2kd+2k-4)} \qquad (t=1,2),$$

$$\mu = \frac{\prod_{t=1}^2(k^t(d-2)+2)}{q(nd-2kd+2k-4)};$$

compare the proof of part (i).

Using (2.21), we obtain then

$$P_{k^1,k^2} = \sum_{D^1,D^2} P(2D^1, 2D^2)$$

$$\sim \frac{d^k((n-k)d)! 2^{kd-k+2} \prod_{s=1}^2 (k^s(d-1))! p^{k-1} q^{kd-2k+4}}{2^{nd/2}(nd-1)!! \prod_{t=1}^2 (k^t(d-2)+2)! \cdot (nd/2-kd+k-2)!}$$

$$\times \exp\left[\frac{(k(d-2))^2}{2qnd}\right].$$

Here

$$\frac{((n-k)d)!}{2^{nd/2}(nd-1)!!(nd/2-kd+k-2)!}$$

$$= 2^{-(kd-k+2)}(nd)^{-k+2} \exp\left[\frac{k^2}{n} - \frac{k^2}{2nd} - \frac{(k(d-1))^2}{2qnd} + O(k/n+k^3/n^2)\right];$$

compare part (i). Combining the last two equations and

$$\binom{n}{k^1, k^2, n-k} = \frac{n^k}{k^1! k^2!} \exp\left(-\frac{k^2}{2n} + O(k/n+k^3/n^2)\right),$$

we obtain

$$E_{k^1,k^2} = n\, k^1, k^2, n - k P_{k^1,k^2}$$

$$\sim \prod_{t=1}^2 \left[ndp^{k^t-1} q^{k^t d - 2(k^t-1)} \frac{(k^t(d-1))!}{k^t!(k^t d - 2(k^t-1))!}\right] e^{\alpha k^2/n}.$$

This relation and (2.3), (2.4) imply (2.6). □

COROLLARY 1. *Suppose that $p - p^* \to 0$. Let $\omega = \omega(n) \to \infty$, such that*
$$\omega = o(n^{1/3}), \qquad \omega \geq n^{1/3}|p - p^*|.$$

*Then*
$$\lim_{n \to \infty} P(L_n^{(1)} \geq n^{2/3}/\omega^2) = 1.$$

For
$$n^{1/3}|p - p^*| \leq (\log n)^a,$$
the corollary implies that whp $L_n^{(1)}$ is at least $n^{2/3}(\log n)^{-2a}$.

PROOF OF COROLLARY 1. Let $k_n = \lceil n^{2/3}/\omega^2 \rceil$. Introduce $X_n$ the total number of tree components of $MG_{np}$ of size $k \in [k_n, 2k_n]$. Then, by Lemma 1(i),

$$(2.22) \quad \mathrm{E}[X_n] = \sum_{k=k_n}^{2k_n} E_k \sim nd\frac{q^2}{p} \sum_{k=k_n}^{2k_n} (pq^{d-2})^k \frac{(k(d-1))!}{k!(kd - 2(k-1))!} e^{\alpha k^2/n}.$$

Here, by (2.5) and the condition of Corollary 1,
$$(2.23) \qquad \alpha \frac{k^2}{n} = O(\omega^{-3}).$$

Furthermore, by the Stirling formula for factorials,
$$\frac{(k(d-1))!}{k!(kd - 2(k-1))!} = \frac{1}{k^{5/2}} \left[\frac{(d-1)^{d-1}}{(d-2)^{d-2}}\right]^k \exp(O(1)),$$
if $k \to \infty$. Observe that
$$\frac{(d-1)^{d-1}}{(d-2)^{d-2}} = [p^*(q^*)^{d-2}]^{-1},$$
and that the function
$$f(p) = pq^{d-2} = p(1-p)^{d-2}$$
attains its absolute maximum at $p = p^*$, with $f''(p^*) > 0$. Hence
$$\frac{f^k(p)}{f^k(p^*)} = \exp[O(k(p - p^*)^2)] = \exp(O(1)),$$
uniformly for $k \in [k_n, 2k_n]$. Therefore $\mathrm{E}[X_n]$, given by (2.22), is of an exact order
$$n \sum_{k=k_n}^{2k_n} k^{-5/2} \sim cnk_n^{-3/2} \sim c\omega^3 \to \infty.$$

Further, by Lemma 1(ii) and (2.23),

$$\mathrm{E}[X_n(X_n - 1)] = \sum_{k_n \leq k^1, k^2 \leq 2k_n} E_{k^1, k^2} \sim \left(\sum_{k=k_n}^{2k_n} E_k\right)^2.$$

Therefore

$$\mathrm{var}[X_n] = o((\mathrm{E}[X_n])^2),$$

and, by Chebyshev's inequality,

$$\lim_{n \to \infty} \mathrm{P}\left(\left|\frac{X_n}{\mathrm{E}[X_n]} - 1\right| \leq \varepsilon\right) = 1.$$

Thus $\mathrm{P}(X_n > 0) \to 1$, and it remains to notice that $L_n^{(1)} \geq n^{2/3}/\omega^2$ if $X_n > 0$.
□

**3. Percolation process on the random pairing.** The upper bound for $L_n^{(1)}$ in Theorem 2, and Theorems 1, 3 will be proved via analysis of a stochastic process which describes growth of a component of $MG_{np}$ containing a given vertex, or more generally, growth of a set of vertices that can be reached from a given set of vertices via open edges. This growth process for $MG_{np}$ is defined through its counterpart on the random pairing $\mathcal{P}_n$; compare [4]. It will be convenient to assign the labels $1, \ldots, nd$ to the points of $\bigcup_{i \in [n]} S_i$.

3.1. *Percolation as a Markov chain.* Let $\mathcal{A}(0) \subset [n]$ be given. We interpret $\mathcal{A}(0)$ as an initial set of active vertices. Then $\mathcal{I}(0) = [n] \setminus \mathcal{A}(0)$ is an initial set of inactive vertices. For each $i \in \mathcal{A}(0)$ [$i \in \mathcal{I}(0)$, resp.], we have the set $S_i(0)$, of cardinality $d$, of active (inactive, resp.) points. We are about to define a process $\{\mathcal{A}(t), \mathcal{I}(t), S_i(t), i \in \mathcal{A}(t) \cup \mathcal{I}(t)\}_{t \geq 0}$, where $S_i(t) \subseteq S_i$. Naturally, $\bigcup_{i \in \mathcal{A}(t)} S_i(t)$, $\bigcup_{i \in \mathcal{I}(t)} S_i(t)$ are called the set of currently active points, and the set of currently inactive points, after $t$ steps.

At step $t+1$ we (a) choose an active point $s' \in \bigcup_{i \in \mathcal{A}(t)} S_i(t)$, with the smallest label, say; (b) identify $s''$, the partner of $s'$ in the random pairing $\mathcal{P}_n$ and delete both $s'$ and $s''$ from the two sets $S_{i'}(t)$, $S_{i''}(t)$ that contain them; (c) if $i'' \in \mathcal{I}(t)$ and the pair $(s', s'')$ is open, then vertex $i''$ is pulled from $\mathcal{I}(t)$, ($\mathcal{I}(t+1) = \mathcal{I}(t) \setminus \{i''\}$), and added to $\mathcal{A}(t)$, ($\mathcal{A}(t+1) = \mathcal{A}(t) \cup \{s''\}$), so that all the points in $S_{i''}(t+1) := S_{i''}(t) \setminus \{s''\}$ become active.

Clearly $\mathcal{A}(t)$ is a set of some of the vertices that can be reached in $MG_{np}$ from the set $\mathcal{A}(0)$ by a path of length $t$ or less. It is possible that an inactive set $S_i(t)$ becomes empty, which means that in $MG_{np}$ there are no loops at the corresponding vertex $i$, and all $d$ edges incident to $i$ are closed.

Given the information on the active points and their partners chosen and deleted in the first $t$ steps, the random pairing on the set of $nd - 2t$ remaining points remains uniformly distributed. In other words, $s''$ is chosen uniformly at random in the set $\bigcup_{i\in[n]} S_i(t) \setminus \{s'\}$.

Introduce
$$I_j(t) = |\{i \in \mathcal{I}(t) : |S_i(t)| = j\}|, \qquad 0 \le j \le d,$$

the total number of currently inactive vertices $i$ such that the points from $S_i = S_i(0)$ have been chosen $d - j$ times as partners $s''$ of active points $s'$ during the first $t$ steps. [The reason that those $i$'s are still inactive after $t$ steps is that the corresponding $d - j$ pairs $(s', s'')$ have all been closed.] Introduce also
$$I(t) = \sum_j j I_j(t),$$

the total count of currently inactive points, and
$$A(t) = \sum_{i \in \mathcal{A}(t)} |S_i(t)|,$$

the total number of currently present active points; in particular, $A(0) = d|\mathcal{A}(0)|$. Let $\mathbf{I}(t) = \{I_j(t)\}_{0 \le j \le d}$. It is easy to see that $\{\mathbf{X}(t)\} = \{A(t), \mathbf{I}(t)\}$ is a Markov chain. Assuming $A(t) > 0$, let us compute the one-step transition probabilities. There are three kinds of transitions:

(a) $s''$ is currently active. The (conditional) probability of this transition is
$$\frac{A(t) - 1}{A(t) - 1 + I(t)}.$$

The next state is
$$A(t+1) = A(t) - 2, \qquad \mathbf{I}(t+1) = \mathbf{I}(t).$$

(b) $s''$ is currently inactive, that is, $s'' \in S_i(t)$, $i \in \mathcal{I}(t)$, $|S_i(t)| = j$, for some $j \in [1, d]$. There are two alternatives.

(b1) $(s', s'')$ is open. The probability of this transition is
$$\frac{j I_j(t) p}{A(t) - 1 + I(t)}.$$

The next state is
$$A(t+1) = A(t) + j - 2,$$
$$I_\ell(t+1) = I_\ell(t), \qquad \ell \ne j,$$
$$I_j(t+1) = I_j(t) - 1.$$

(b2) $(s', s'')$ is closed. The probability of this transition is
$$\frac{jI_j(t)q}{A(t) - 1 + I(t)}.$$
The next state is
$$A(t+1) = A(t) - 1, \qquad I_\ell(t+1) = I_\ell(t), \qquad \ell \neq j-1, j,$$
$$I_{j-1}(t+1) = I_{j-1}(t) + 1, \qquad I_j(t+1) = I_j(t) - 1.$$
In all cases,
$$(3.1) \qquad A(t+1) + I(t+1) = A(t) + I(t) - 2,$$
as it should be, of course. If $A(t) = 0$, then $\mathbf{X}(t+1) = \mathbf{X}(t)$, that is, every $(0, \mathbf{I})$ is an absorbing state.

NOTE. This Markov chain is a close relative of the Markov chain for the bootstrap site percolation in [4]. There are substantial differences though. One is that in [4] $|\mathcal{A}(0)|$, the number of active vertices at the start, was of order $n$, while here we will have to consider $|\mathcal{A}(0)| = 1$ or $2$, that is, to start with just one or two active vertices. With $|\mathcal{A}(0)|$ that small, the process runs out of active vertices relatively soon, with probability $P_n$ bounded away from zero, or even approaching 1, if $n^{1/3}(p - p^*) \not\to \infty$. Since there are $n$ vertices, our task basically is to find a way to handle something like $n(1 - P_n)$ for various ranges of $p$.

To continue, we average over the three possibilities, (a), (b1) and (b2), and obtain the equations for the conditional expectations $\mathbb{E}[\mathbf{X}(t+1)|\mathbf{X}(t)] = \mathbb{E}[\mathbf{X}(t+1)|\circ]$: if $A(t) > 0$, then
$$\mathbb{E}[A(t+1)|\circ] = A(t) + \frac{A(t) - 1}{A(t) - 1 + I(t)}(-2)$$
$$+ \sum_j \frac{jI_j(t)p}{A(t) - 1 + I(t)}(j - 2)$$
$$(3.2) \qquad + \sum_j \frac{jI_j(t)q}{A(t) - 1 + I(t)}(-1),$$
$$\mathbb{E}[I_j(t+1)|\circ] = I_j(t) + \frac{jI_j(t)}{A(t) - 1 + I(t)}(-1)$$
$$+ \frac{(j+1)I_{j+1}(t)q}{A(t) - 1 + I(t)}, \qquad 0 \leq j \leq d,$$
$I_{d+1}(t) \equiv 0$. [$I_j(t+1) = I_j(t) - 1$ when $s''$ belongs to an inactive set of cardinality $j$, no matter whether $(s', s'')$ is open or closed.] We will denote the right-hand side of (3.2) by $\mathbf{X}(t) + \mathbf{R}(\mathbf{X}(t))$.

### 3.2. *Differential equations approximation.*

3.2(a). As we remarked earlier (note following Lemma 1), whp the number of isolated vertices in $MG_{np}$ is of order $n$ exactly. This means that whp $I(t)$ is of order $n$ for every $t$, suggesting to scale $\mathbf{X}(t)$ by $n$. Another argument in favor of such scaling is that, if not for $-1$ in $A(t) - 1$, $R(\mathbf{X})$ would have been homogeneous, zero degree, vector-function of $\mathbf{X}$. Besides

$$1 = (t+1) - t = n\left(\frac{t+1}{n} - \frac{t}{n}\right).$$

Hoping that, when it matters, $\mathbf{X}(t+1)$ is relatively close to its conditional expectation $\mathbb{E}[\mathbf{X}(t+1)|\circ]$, (3.2) leads us to conjecture that the random sequence $\{\mathbf{X}(t)\}$ is well approximated by $\{n\mathbf{x}(t/n)\}$. Here $\mathbf{x}(\tau) = (a(\tau), \mathbf{i}(\tau))$ is the solution of the (deterministic) differential equations

$$(3.3) \quad a' = \frac{a}{a+i}(-2) + p\sum_j \frac{ji_j}{a+i}(j-2) + q\sum_j \frac{ji_j}{a+i}(-1),$$

$$i'_j = \frac{q(j+1)i_{j+1} - ji_j}{a+i}, \qquad 0 \le j \le d,$$

where $i_{d+1}(\tau) \equiv 0$, subject to initial conditions

$$(3.4) \quad a(0) = \frac{A(0)}{n}, \qquad i_j(0) = 0,\ 0 \le j \le d-1, \qquad i_d(0) = 1 - \frac{A(0)}{nd}.$$

We will denote the right-hand side of (3.3) by $\mathbf{R}_0(\mathbf{x})$.

Let us solve (3.3). It follows from the equations that

$$(a(\tau) + i(\tau))' = -2;$$

in view of (3.1) this is hardly surprising. As $a(0) + i(0) = d$, we have

$$a(\tau) + i(\tau) = d - 2\tau.$$

Therefore we consider $\tau < d/2$ only. Introducing $u = \ln(a(\tau) + i(\tau))^{-1/2}$, and $f_j(u) := i_j(\tau)$, we obtain a system of *linear* (birth-and-death type) equations

$$\frac{df_j}{du} = q(j+1)f_{j+1} - jf_j, \qquad 0 \le j \le d,\ f_{d+1}(u) \equiv 0.$$

By a backward induction on $j$ it follows that, for all $u_1$, $u_2$,

$$(3.5) \quad f_j(u_2) = e^{-ju}\sum_{r=j}^{d} q^{r-j}\binom{r}{j}(1-e^{-u})^{r-j}f_r(u_1), \qquad u = u_2 - u_1.$$

Or, setting $\mathbf{f}(\cdot) = (f_0(\cdot), \ldots, f_d(\cdot))$,

$$(3.6) \quad \mathbf{f}(u_2) = M(u)\mathbf{f}(u_1), \qquad M_{jr}(u) = e^{-ju}q^{r-j}\binom{r}{j}(1-e^{-u})^{r-j}\mathbf{1}_{\{r \ge j\}}.$$

It follows that

(3.7) $$M(v_1 + v_2) = M(v_1)M(v_2), \qquad M(-v) = M(v)^{-1}.$$

For $u_2 = \ln(a(\tau) + i(\tau))^{-1/2}$, $u_1 = \ln(a(0) + i(0))^{-1/2}$, (3.5) yields

(3.8)
$$i_j(\tau) = y(\tau)^j \sum_{r=j}^{d} q^{r-j} r j (1 - y(\tau))^{r-j} i_r(0), \qquad 0 \leq j \leq d,$$

$$y(\tau) := \left(\frac{a(\tau) + i(\tau)}{a(0) + i(0)}\right)^{1/2} = \left(1 - \frac{2\tau}{a(0) + i(0)}\right)^{1/2}.$$

For the initial conditions (3.4), the formula (3.8) becomes

(3.9)
$$i_j(\tau) = dj y(\tau)^j [q(1 - y(\tau))]^{d-j} i_d(0), \qquad 0 \leq j \leq d,$$

$$y(\tau) = \left(1 - \frac{2\tau}{d}\right)^{1/2}.$$

Therefore

(3.10) $$i(\tau) = \sum_j j i_j(\tau) = dy(\tau)(py(\tau) + q)^{d-1} i_d(0),$$

and

(3.11)
$$a(\tau) = (a(\tau) + i(\tau)) - i(\tau) = S(y(\tau), p),$$

$$S(y(\tau), p) := dy(\tau)[y(\tau) - (py(\tau) + q)^{d-1} i_d(0)].$$

We will need to consider $A(0)$ equal to $d$ or $2d$ only, corresponding to a single active vertex, or to two active vertices at $t = 0$. So in (3.10)–(3.11), $i_d(0) = 1 - n^{-1}$ or $i_d(0) = 1 - 2n^{-1}$. $S(y, p) = 0$ has two roots, 0 and $\hat{y} = \hat{y}(p) \in (0, 1)$. Intuitively we anticipate that the random process will run out of active points, thus will terminate, at a moment close to $n\hat{\tau}$,

(3.12) $$\hat{\tau} = \hat{\tau}(p) = \frac{d}{2}(1 - \hat{y}^2);$$

see (3.11), second line. $\hat{\tau} = O(1/n)$ for $p < p^*$, and $\hat{\tau}$ is bounded away from zero if $p > p^*$. Our task is to find a rigorous argument that will also cover the case $p - p^* \to 0$.

3.2(b). For $u_2 = 0$, $u_1 = \ln(a(\tau) + i(\tau))^{-1/2}$, (3.5) becomes

$$F_j(\mathbf{x}(\tau)) = i_j(0),$$

(3.13) $$F_j(\mathbf{x}) := y^{-j} \sum_{r=j}^{d} q^{r-j} \binom{r}{j} (1 - y^{-1})^{r-j} i_r, \qquad 0 \leq j \leq d,$$

$$y = (a + i)^{1/2}.$$

In a vector-matrix form, the second line in (3.13) can be expressed as

(3.14) $$\mathbf{F}(\mathbf{x}) = M(\ln(a+i)^{1/2})\mathbf{i}, \qquad \mathbf{x} = (a, \mathbf{i}).$$

While (3.8) solves the initial value problem, (3.13) provides a collection of $d+1$ functions of $\mathbf{x}$, which are the integrals of (3.3). That is, we have

(3.15) $$\mathbf{grad}^* F_j(\mathbf{x}) \cdot \mathbf{R}_0(\mathbf{x}) \equiv 0.$$

Next we use the integrals $F_j(\mathbf{x})$ and (3.15) to construct a collection of exponential supermartingales, which will allow us to bound probability of "large" deviations of $\mathbf{X}(t)$ from $n\mathbf{x}(t/n)$.

LEMMA 2. *Let $\gamma = \gamma(n) \to \infty$ and $\gamma = O(n)$. Introduce*
$$Q_j(t) = \exp\{\gamma[F_j(n^{-1}\mathbf{X}(t)) - F_j(n^{-1}\mathbf{X}(0))]\}, \qquad 0 \le j \le d.$$
*There exists an absolute constant $c > 0$ such that*
$$\mathbb{E}[Q_j(t+1)|\circ] \le (1 + c(\gamma/n)^2)Q_j(t),$$
*whenever $I(t) \ge nq^d/2$.*

NOTE. $I(t)$ is bounded below by the number of isolated vertices in $MG_{np}$, which is asymptotic, in probability, to $nq^d$. So whp the condition in the lemma holds for all $t$.

PROOF OF LEMMA 2. Obviously we need to consider the case $A(t) > 0$ only. $F_j(\mathbf{x})$ is certainly twice continuously differentiable for $a + i \ge q^d/2$. Therefore, since $\mathbf{X}(t+1) - \mathbf{X}(t)$ is uniformly bounded,

$$\begin{aligned}\frac{Q_j(t+1)}{Q_j(t)} &= \exp\{\gamma[F_j(n^{-1}X(t+1)) - F_j(n^{-1}X(t))]\}\\ &= \exp\{\gamma[(\mathbf{grad}F_j(n^{-1}\mathbf{X}(t)))^*(n^{-1}\mathbf{X}(t+1) - n^{-1}\mathbf{X}(t))]\\ &\qquad\qquad + O(\gamma n^{-2})\}\\ &= 1 + \gamma\mathbf{grad}F_j(n^{-1}\mathbf{X}(t))^*(n^{-1}\mathbf{X}(t+1) - n^{-1}\mathbf{X}(t)) + O(\gamma^2/n^2).\end{aligned}$$

Notice that

$$\begin{aligned}\mathbb{E}[(\mathbf{grad}&F_j(n^{-1}\mathbf{X}(t)))^*(n^{-1}\mathbf{X}(t+1) - n^{-1}\mathbf{X}(t))|\circ]\\ &= (\mathbf{grad}F_j(n^{-1}\mathbf{X}(t)))^*\mathbb{E}[n^{-1}\mathbf{X}(t+1) - n^{-1}\mathbf{X}(t)|\circ]\\ &= (\mathbf{grad}F_j(n^{-1}\mathbf{X}(t)))^* n^{-1}\mathbf{R}(\mathbf{X}(t))\\ &= (\mathbf{grad}F_j(n^{-1}\mathbf{X}(t)))^* n^{-1}\mathbf{R}_0(n^{-1}\mathbf{X}(t)) + O(n^{-2})\\ &= O(n^{-2});\end{aligned}$$

see (3.15). Consequently
$$\mathbb{E}[Q_j(t+1)|\circ] = Q_j(t)[1 + O(\gamma/n^2 + \gamma^2/n^2)] \le Q_j(t)(1 + c\gamma^2/n^2). \quad \square$$

LEMMA 3. *Let $T$ be the first time $t$ when either $A(t) = 0$ or $I(t) < nq^d/2$.*
(i) *Then for any fixed $s$, and $z > 0$,*

$$(3.16) \quad P\left\{\max_{t \leq T_s} |F_j(n^{-1}\mathbf{X}(t)) - F_j(n^{-1}\mathbf{X}(0))| > z\right\} \leq 2e^{-\gamma z + cs\gamma^2/n^2},$$

$T_s = T \wedge s$. (ii) *Consequently, for some absolute constant $c_1 > 0$,*

$$(3.17) \quad P\left\{\max_{t \leq T_s} |n^{-1}\mathbf{X}(t) - \mathbf{x}(n^{-1}t)| > z\right\} \leq 2(d+1)e^{-c_1 \gamma z + cs\gamma^2/n^2}.$$

NOTE. For (3.17) to be of any help, $\gamma z$ should substantially exceed $s\gamma^2/n^2$, and moreover $\gamma z$ will have to grow logarithmically with $n$. On the other hand $z$, needs to be much smaller than $x(n^{-1}t)$, for $t$'s in question, otherwise (3.17) is pointless. Selecting $\gamma$, $z$, which meet these conflicting requirements, will depend on the range of $p$ under consideration.

PROOF OF LEMMA 3. (i) Define

$$\hat{Q}_j(t) = \frac{Q_j(t \wedge T_s)}{(1 + c\gamma^2 n^{-2})^{t \wedge T_s}}.$$

By Lemma 2, $\{\hat{Q}_j(t)\}$ is a supermartingale, that is,

$$\mathbb{E}[\hat{Q}_j(t+1)|\circ] \leq \hat{Q}_j(t), \qquad t \geq 0.$$

Let $\mathcal{T}_z$ be the first time $t \leq T_s$ such that

$$F_j(n^{-1}\mathbf{X}(t)) - F_j(n^{-1}\mathbf{X}(0)) > z,$$

and set $\mathcal{T}_z = T_s + 1$ if no such $t$ exists. Now $\mathcal{T}_z$ is a stopping time. So, using the Optional Sampling Theorem (Durrett [14], Chapter 4),

$$\mathbb{E}[\hat{Q}_j(\mathcal{T}_z)] \leq \hat{Q}_j(0) = 1.$$

On the event $\{\mathcal{T}_z \leq T_s\}$ we have

$$\hat{Q}_j(\mathcal{T}_z) \geq \frac{e^{\gamma z}}{(1 + c\gamma^2 n^{-2})^s} \geq e^{\gamma z - cs\gamma^2 n^{-2}}.$$

So

$$P\left\{\max_{t \leq T_s}(F_j(n^{-1}\mathbf{X}(t)) - F_j(n^{-1}\mathbf{X}(0))) > z\right\}$$
$$\leq P\{\hat{Q}_j(\mathcal{T}_z) \geq e^{\gamma z - cs\gamma^2 n^{-2}}\}$$
$$\leq \frac{\mathbb{E}[\hat{Q}_j(\mathcal{T}_z)]}{e^{-\gamma z + cs\gamma^2 n^{-2}}} \leq e^{\gamma z - cs\gamma^2 n^{-2}}.$$

In exactly the same way we obtain

$$\mathrm{P}\left\{\max_{t\leq T_s}(-F_j(n^{-1}\mathbf{X}(t))+F_j(n^{-1}\mathbf{X}(0)))>z\right\}\leq e^{-\gamma z+cs\gamma^2 n^{-2}}.$$

So (3.16) follows.

(ii) According to the definition of $F_j$'s, we have thus proved that

$$\mathrm{P}\left\{\max_{t\leq T_s}\|M(U(t))(n^{-1}\mathbf{I}(t))-M(U(0))(n^{-1}\mathbf{I}(0))\|\leq z\right\}$$
(3.18)
$$\geq 1-(d+1)e^{-\gamma z+cs\gamma^2 n^{-2}}$$

($\|\mathbf{f}\|:=\max_j|f_j|$). Here $U(t)=\ln(n^{-1}(A(t)+I(t)))^{1/2}$ and $M(\cdot)$ is given by (3.6). On the event $E_n$ in (3.18) we have: for $t\leq T_s$,

$$M(U(t))(n^{-1}\mathbf{I}(t))-M(U(0))(n^{-1}\mathbf{I}(0))=\mathbf{Z},\qquad \|\mathbf{Z}\|\leq z.$$

Now $M(U(t))^{-1}=M(-U(t))$, and

$$U(t)\leq \ln(n^{-1}(A(0)+I(0)))^{1/2}=\ln d^{1/2}<\infty.$$

Using the definition of $M(\cdot)$ in (3.6), we see then that $\|M(-U(t))\|\leq m$, for some absolute constant $m>0$. So on $E_n$,

$$mz\geq \|n^{-1}\mathbf{I}(t)-M^{-1}(U(t))M(U(0))n^{-1}\mathbf{I}(0)\|$$
$$=\|n^{-1}\mathbf{I}(t)-M(U(0)-U(t))n^{-1}\mathbf{I}(0)\|.$$

Here

$$U(0)-U(t)=\ln\left(\frac{n^{-1}A(0)+n^{-1}I(0)}{n^{-1}A(t)+n^{-1}I(t)}\right)^{1/2}=\ln\left(1-\frac{2t}{nd}\right)^{-1/2},$$

so

$$M(U(0)-U(t))n^{-1}\mathbf{I}(0)=\mathbf{i}(t/n);$$

see (3.8), (3.9). Therefore it follows from (3.18) that

$$\mathrm{P}\left\{\max_{t\leq T_s}\|n^{-1}\mathbf{I}(t)-\mathbf{i}(t/n)\|\leq mz\right\}$$
$$\geq 1-(d+1)e^{-\gamma z+cs\gamma^2 n^{-2}}.$$

This together with

$$n^{-1}A(t)+n^{-1}I(t)=a(t/n)+i(t/n)=d-2t/n$$

imply (3.17). □

**4. Proof of Theorem 1 (Subcritical case).** Let $C_n$ be size of the component of $MG_{np}$ which contains a given vertex, 1 say. This component can be determined through the percolation process with the starting active set $\mathcal{A}(0) = \{1\}$. In Lemma 3 we introduced the stopping time $T$, the first $t$ such that either $A(t)$ (the number of currently active points) is zero, or $I(t)$ (the number of currently inactive points) is below $nq^d/2$. Clearly $C_n \leq 1+T$. Set in Lemma 3

$$\gamma = \alpha(p^* - p)n, \qquad s = \beta(p^* - p)^{-2}\ln n,$$

$\alpha > 0$, $\beta > 0$ to be specified shortly. Pick $z > 0$ such that

(4.1) $$\gamma z \geq 2cs\frac{\gamma^2}{n^2} \iff z \geq 2c\alpha\frac{s(p^* - p)}{n}$$

in which case

$$\gamma z - cs\frac{\gamma^2}{n^2} \geq c\alpha^2\beta \ln n.$$

Then, by Lemma 3,

(4.2) $$\mathrm{P}\left\{\max_{t \leq T_s}|n^{-1}\mathbf{X}(t) - \mathbf{x}(n^{-1}t)| > z\right\} = O(n^{-c\alpha^2\beta}) = o(n^{-1}),$$

if $c\alpha^2\beta > 1$, the restriction we meet by selecting $\beta > c^{-1}\alpha^{-2}$. On the event $E_{n1}$, complementary to that in (4.2), we have

(4.3) $$|n^{-1}A(T_s) - a(T_s/n)| \leq z.$$

Here, by (3.11),

(4.4) $$a(\tau) = dy(\tau)[y(\tau) - (py(\tau) + q)^{d-1}(1 - n^{-1})],$$
$$y(\tau) = (1 - 2\tau/d)^{1/2}.$$

Suppose that, on $E_{n1}$, $T > s$. Then (4.3) implies that $z + a(s/n)$ must be positive. Since

(4.5) $$\frac{s}{n} = \beta\frac{\ln n}{n(p^* - p)^2} = \beta\frac{\ln n}{n^{1/3}(n^{1/3}(p^* - p))^2} \to 0,$$

as $n^{1/3}(p^* - p) \to \infty$, we have

$$1 - y(s/n) = 1 - (1 - 2s/(nd))^{1/2} \sim s/(nd).$$

Therefore, by convexity of $y - (py + q)^{d-1}(1 - n^{-1})$,

$$a(s/n) \leq dy(s/n)\{n^{-1} + [1 - p(d - 1)(1 - n^{-1})](y(s/n) - 1)\}$$
$$= -dy(s/n)[(p^* - p)(1 - y(s/n)) + O(n^{-1})]$$
$$\sim -\frac{s(p^* - p)}{n},$$

as $(p^* - p)s \to \infty$. So $z + a(s/n) < 0$ provided that

(4.6) $$z \leq \frac{(p*-p)s}{2n}.$$

We meet the restrictions (4.1) and (4.6) by choosing, for instance,

$$z = 0.5(2c\alpha + 0.5)s(p^* - p).$$

With $z$ so defined we come to the conclusion that $T \leq s$ on the event $E_{n1}$. Thus

$$P(C_n > 1 + s) \leq P(T > s) \leq P(E_{n1}^c) = o(n^{-1}).$$

Using the union bound we see that whp $MG_{np}$ does not have a component of size exceeding $1 + \beta(p^* - p)^{-2} \log n$, that is,

(4.7) $$L_n^{(1)} \leq 1 + \beta(p^* - p)^{-2} \ln n.$$

This completes the proof of Theorem 1.

NOTE. Observe that $s/n \to 0$ if $n^{1/2}(p^* - p) \gg \ln n$; see (4.5). That is, the statement holds under this weaker condition. However, Theorem 2 makes it clear that when $p^* - p$ is of order $n^{-a}$, $a \in (1/3, 1/2)$, the upper bound for $L_n^{(1)}$ is roughly $n^{2/3}$, thus smaller than our current $s$ by $n^{2(a-1/3)}$ factor.

**5. Proof of Theorem 2 (Nearcritical case). Upper bound.** The argument parallels the previous proof. This time we set

(5.1) $$\gamma = n^{2/3}(\ln n)^\alpha, \qquad s = n^{2/3}(\ln n)^\beta,$$

for $\alpha > 0$, $\beta > 0$ to be specified. The counterpart of (4.1) is the condition

(5.2) $$z \geq 2c\frac{s\gamma}{n^2} \quad \Longleftrightarrow \quad z \geq 2c\frac{(\ln n)^{\alpha+\beta}}{n^{2/3}}.$$

From (4.4) it follows that

$$a(s/n) = dy(s/n)\left[n^{-1} + (1 - p(d-1)(1 - n^{-1}))(y(s/n) - 1)\right.$$
$$\left. - p\binom{d-1}{2}(1 - n^{-1})(y(s/n) - 1)^2 + O((1 - y(s/n))^3)\right]$$
$$\sim -dp\,d - 12(y(s/n) - 1)^2 \sim -d^{-1}p\binom{d-1}{2}\frac{s^2}{n^2},$$

provided that $s/n \gg |p^* - p|$. Since $|p^* - p| \leq n^{-1/3}(\ln n)^a$, we satisfy this condition by choosing in (5.1) $\beta > a$. Then $z + a(s/n) < 0$ if

(5.3) $$z \leq (d-1)(d-2)\frac{s^2}{n^2} = (d-1)(d-2)\frac{(\ln n)^{2\beta}}{n^{2/3}}.$$

Comparing (5.2) and (5.3) we see that

$$z = 2c\frac{(\ln n)^{\alpha+\beta}}{n^{2/3}}$$

satisfies both conditions if $\alpha < \beta$. For this choice of $z$,

$$\mathrm{P}(E^c_{n1}) \le e^{-\gamma z/2} = e^{-c(\ln n)^{2\alpha+\beta}} = o(n^{-1}),$$

if $2\alpha + \beta > 1$. If $\beta > \max\{1/3, a\}$, we make $2\alpha + \beta > 1$ by selecting $\alpha$ sufficiently close to $\beta$ from below. For these $\alpha$, $\beta$, we have thus proved that

$$\lim_{n \to \infty} \mathrm{P}(L^{(1)}_n \ge 1 + n^{2/3}(\ln n)^\beta) = 0.$$

**6. Proof of Theorem 3 (Supercritical case).** Now

(6.1) $$\omega = \omega(n) := n^{1/3}(p - p^*) \ge (\ln n)^{a(n)},$$

where $a(n) \to \infty$ however slowly.

*Preliminaries.* We will have to consider two cases, with one and two initially active vertices. So this time

(6.2)
$$a(\tau) = dy(\tau)[y(\tau) - (py(\tau) + q)^{d-1}(1 - m/n)],$$
$$y(\tau) = (1 - 2\tau/d)^{1/2},$$

where $m = 1$, or $m = 2$.

Let $\hat{y} = \hat{y}(p)$ be the root of $f_{m/n}(y) = 0$,

$$f_{m/n}(y) := y - (py + q)^{d-1}(1 - m/n),$$

so that

$$\hat{\tau} = \hat{\tau}(p) = \frac{d}{2}(1 - \hat{y}^2)$$

is the root of $a(\tau) = 0$ in $(0, d/2)$. Since $f_{m/n}(y) \ge f_0(y)$ and

(6.3) $$-(d-1)(d-2)p^2 \le \frac{d^2 f_{m/n}(y)}{dy^2} \le -\frac{1}{2}(d-1)(d-2)p^2 q^d,$$

it follows easily that

(6.4) $$|\hat{y} - \pi| = O(n^{-1/2}),$$

where $\pi \in (0, 1)$ is the root of $f_0(y) = 0$ in $(0, 1)$. From (6.3) for $m = 0$, it follows directly that $1 - \pi$ is of order $p - p^*$ exactly, which we express by writing

$$1 - \pi \asymp p - p^*.$$

Then, by (6.4),

(6.5) $$1 - \hat{y} \asymp p - p^*, \qquad \hat{\tau} \asymp p - p^*,$$

too, as $p - p^* \gg n^{-1/3}$. We will also need bounds for $a(\tau)$ when $\tau$ is (relatively) close to 0 or $\hat{\tau}$. The convex function $f_{m/n}(y)$ attains its maximum at a point $\tilde{y} \in (\hat{y}, 1)$, and a simple calculus shows that

$$\frac{1-\tilde{y}}{1-\hat{y}} = \frac{1}{2} + O(p - p^*) \implies \tilde{y} - \hat{y} \asymp p - p^*.$$

Using this fact, $f'_{m/n}(\tilde{y}) = 0$, and (6.3), we obtain that

$$f'_{m/n}(y) \asymp \begin{cases} (p - p^*), & \dfrac{|y - \hat{y}|}{p - p^*} \leq \varepsilon_0, \\ -(p - p^*), & \dfrac{1 - y}{p - p^*} \leq \varepsilon_0, \end{cases}$$

if $\varepsilon_0 > 0$ is sufficiently small. Therefore, as $f_{m/n}(\hat{y}) = f_{m/n}(1) = 0$,

$$f_{m/n}(y) \asymp -(y - \hat{y})(p - p^*), \qquad \frac{|y - \hat{y}|}{p - p^*} \leq \varepsilon_0,$$
$$f_{m/n}(y) \asymp (1 - y)(p - p^*), \qquad \frac{1 - y}{p - p^*} \leq \varepsilon_0.$$

Since $a(t) = dy(\tau) f_{m/n}(y(\tau))$, a little reflection shows then that

(6.6) $$\min\{a(\tau) : \tau \in [\varepsilon(p - p^*), \hat{\tau} - \varepsilon(p - p^*)]\} \geq c\varepsilon(p - p^*)^2,$$
$$a(\hat{\tau} + \varepsilon(p - p^*)) \leq -c\varepsilon(p - p^*)^2,$$

if $\varepsilon > 0$, $c > 0$ are sufficiently small.

*Running time dichotomy.* Let $m = 1$. In Lemma 3 set

(6.7) $$\gamma = n^{2/3}, \quad s = n\hat{\tau}(1 + \beta\omega^{-1}),$$

$\beta > 0$ being fixed. The analogue of the conditions (4.1), (5.2) is

$$\gamma z \geq 2cs\frac{\gamma^2}{n^2} \iff z \geq 2c\hat{\tau}(1 + \beta\omega^{-1})n^{-1/3},$$

which is met (for $n$ large enough) if

(6.8) $$z \geq c'\omega n^{-2/3}$$

[see (6.5)], where $c'$ is independent of $\beta$. Further, using (6.5) and the second line in (6.6),

$$a(s/n) = a(\hat{\tau} + \beta\hat{\tau}\omega^{-1}) \leq -c''\frac{\beta\omega}{n^{2/3}}.$$

So, picking $\beta > 1$ sufficiently large and setting $z = \beta^{1/2}\omega n^{-2/3}$, say, we ensure that $z$ satisfies (6.8) and

(6.9) $$z + a(s/n) < 0.$$

Since $z = z(\beta)$ satisfies (6.8), by Lemma 3(ii) we obtain: there exist $c_i = c_i(\beta) > 0$, such that

(6.10) 
$$\begin{aligned}
&P\{\|n^{-1}\mathbf{X}(T_s) - \mathbf{x}(n^{-1}T_s)\| > z\} \\
&\geq 1 - \exp(-c_1\gamma z) \\
&\geq 1 - \exp(-c_2\omega) = 1 - O(n^{-K}) \qquad \forall K > 0.
\end{aligned}$$

Let us denote the event in (6.10) by $B_n$. On $B_n \cap \{T > s\}$

$$z + a(s/n) \geq n^{-1}A(T_s) > 0,$$

which contradicts (6.9). Hence

$$B_n \in \{T \leq s\} \cap \{|n^{-1}\mathbf{X}(T) - \mathbf{x}(n^{-1}T)| \leq z\}.$$

In particular, since $A(T) = 0$,

(6.11) $$z \geq |n^{-1}A(T) - a(n^{-1}T)| = |a(n^{-1}T)|.$$

Suppose that

$$\frac{\beta}{c\omega}(p - p^*) \leq \frac{T}{n} \leq \hat{\tau} - \frac{\beta}{c\omega}(p - p^*).$$

Then, by the first line in (6.6),

$$a(n^{-1}T) \geq \frac{\beta}{\omega}(p - p^*)^2 = \beta\frac{\omega}{n^{2/3}},$$

and, by (6.11) and $z = \beta^{1/2}\omega n^{-2/3}$, we must have $\beta^{1/2} \geq \beta$, which is impossible since $\beta > 1$. Therefore

$$B_n \subseteq \{T \leq s\} \cap \{n^{-1}T \notin [c^{-1}\beta n^{-1/3}, \hat{\tau} - c^{-1}\beta n^{-1/3}]\},$$

as $p - p^* = n^{-1/3}\omega$. Or, recalling that $s = n\hat{\tau}(1 + \beta\omega^{-1})$, $\hat{\tau} \asymp p - p^*$,

(6.12) $B_n \subseteq \{n^{-1}T \in [\hat{\tau} + \alpha n^{-1/3}, \hat{\tau} - \alpha n^{-1/3}]\} \cup \{n^{-1}T \leq \alpha n^{-1/3}\},$

for some constant $\alpha > 0$. In words, with probability superpolynomially close to 1, the random percolation —triggered by a fixed initially active vertex— either stops after $O(n^{2/3})$ steps, or runs much longer, for about $n\hat{\tau} \asymp \omega n^{2/3}$ steps.

*... and uniqueness of a giant component.* On
$$B_n^* := B_n \cap \{n^{-1}T \in [\hat{\tau} + \alpha n^{-1/3}, \hat{\tau} - \alpha n^{-1/3}]\}$$
we have [see (3.9)]
$$n^{-1}\sum_{j=0}^{d} I_j(T) = \sum_{j=0}^{d} i_j(n^{-1}(T)) + O(\omega n^{-2/3})$$
$$= (py(n^{-1}T) + q)^d(1 - n^{-1}) + O(\omega n^{-2/3})$$
$$= (p\pi + q)^d + O(n^{-1/3} + \omega n^{-2/3}).$$

Therefore $C_n$, the size of the component that contains our initially active vertex, on $B_n^*$ satisfies

(6.13)
$$C_n = n - \sum_{j=0}^{d} I_j(T)n = n[1 - (p\pi + q)^d] + O(n^{2/3})$$
$$\asymp \omega n^{2/3} = n(p - p^*).$$

And, of course, $C_n \leq 1 + T \leq 1 + \alpha n^{2/3}$ on $B_n \cap \{T \leq \alpha n^{2/3}\}$. Thus, with probability $1 - O(n^{-K})$, $\forall K > 0$, either $|T - n\hat{\tau}| \leq an^{2/3}$ and $C_n \sim n[1 - (p\pi + q)^d]$, or $T \leq 1 + \alpha n^{2/3}$.

We observe now that no changes whatsoever are needed to show validity of this claim for $\tilde{T}$ and $\tilde{C}_n$, the running time and the terminal number of active vertices of the percolation process started by *two* fixed vertices $u$ and $v$. That is, with probability $1 - O(n^{-K})$, $\forall K > 0$,

(6.14)
$$\tilde{C}_n \leq 1.1n[1 - (p\pi + q)^d].$$

And the event (6.14) excludes a possibility that $u$ and $v$ belong to two different components, each of size relatively close to $n[1 - (p\pi + q)^d]$. Picking $K > 2$, we obtain then: with probability $1 - O(n^{-(K-2)})$, there are no two vertices $u, v \in [n]$ such that the two components containing them are *disjoint*, and both are of size of order $n(p - p^*)$. Putting it differently, whp there may exist at most one "giant" component, and all other components have size $O(n^{2/3})$ at most.

*Other components....* Let us show whp that nongiant components are actually much smaller, of order $n^{2/3}\omega^{-2}(\ln n)^{1+o(1)}$.

To this end, define a sequence $\{\gamma_k, s_k, z_k\}$. For $k \geq 0$,
$$\gamma_k = \omega^{\sigma_k} n^{2/3}, \qquad s_k = \alpha\omega^{-\beta_k} n^{2/3},$$
and

(6.15)
$$z_k \geq 2cs_0 \frac{\gamma_k}{n^2} = 2c\alpha\omega^{\sigma_k - \beta_k},$$

$\sigma_k > 0$, $\beta_k \geq 0$. In particular, $\beta_0 = 0$, so that $s_0 = \alpha n^{2/3}$ is the (probabilistic) upper bound for all nongiant components already proved. Further, $\beta_k$ increases with $k$, and

$$\sigma_k = \frac{\beta_k}{2} + \frac{\rho \ln \ln n}{\ln \omega}, \qquad \rho > 1. \tag{6.16}$$

So

$$\gamma_k z_k - cs_k \frac{\gamma_k^2}{n^2} = c\alpha \omega^{2\sigma_k - \beta_k} = (\ln n)^\rho.$$

Applying Lemma 3 with $\gamma_k$, $s_k$ and $z_k$, we obtain that

$$|n^{-1}\mathbf{X}(t) - \mathbf{x}(n^{-1}t)| \leq z_k \qquad \forall t \leq T_{s_k}(= \min\{T, s_k\}), \tag{6.17}$$

with probability $1 - O(e^{-(\ln n)^\rho})$. Now

$$\frac{n^{-1}s_k}{p - p^*} = O(\omega^{-1}) \to 0;$$

so, using (6.6),

$$\min\{a(t) : t \in [s_{k+1}/n, s_k/n]\} \geq a(n^{-1}s_{k+1}) \asymp \frac{\alpha}{\omega^{\beta_{k+1}-1}}. \tag{6.18}$$

Pick $\rho' \in (1, \rho)$, and set

$$\beta_{k+1} = 1 + \frac{\beta_k}{2} - \frac{\rho' \ln \ln n}{2 \ln \omega}. \tag{6.19}$$

Then [see (6.16)],

$$\beta_{k+1} - 1 - (\beta_k - \sigma_k) = \frac{(\rho - \rho') \ln \ln n}{2 \ln \omega},$$

which means that we can pick $z_k$ satisfying (6.15), and such that

$$z_k = o(a(n^{-1}s_{k+1}));$$

see (6.18). But then $T \notin [s_{k+1}, s_k]$ on the event (6.17), since otherwise, by (6.17), we obtain that

$$a(n^{-1}s_{k+1}) \geq |a(n^{-1}(T))| = |a(n^{-1}T) - n^{-1}A(T)| \leq z_k,$$

contradicting our choice of $z_k$.

The solution of (6.19), subject to $\beta_0 = 0$, is

$$\beta_k = 2 - \frac{\rho' \ln \ln n}{\ln \omega} - 2^{-k+1}\left(1 - \frac{\rho' \ln \ln n}{2 \ln \omega}\right).$$

For $k_0 = \lceil \log_2 \ln \omega \rceil$,

$$\beta_{k_0} = 2 - \frac{(\rho' + o(1)) \ln \ln n}{\ln \omega},$$

so that

$$s_{k_0} = \alpha \omega^{-\beta_{k_0}} n^{2/3} = \alpha \frac{n^{2/3}(\ln n)^{\rho'+o(1)}}{\omega^2},$$

for all $1 < \rho' < \rho$. By the union bound inequality, we have then

$$P\{T \in [s_{k_0}, s_0]\} \leq \sum_{k=0}^{k_0-1} P\{T \in [s_{k+1}, s_k]\}$$
$$= O(e^{-(\ln n)^\rho} k_0) = O(e^{-(\ln n)^\rho} \ln \ln n)$$
$$= O(n^{-K}) \quad \forall K > 0.$$

Now $s_0 = \alpha n^{2/3}$, and we had proved that

$$P\{T \in [\alpha n^{2/3}, n\hat{\tau} - \alpha n^{2/3}]\} = O(n^{-K}).$$

Consequently

$$P\{T \in [\omega^{-2} n^{2/3}(\ln n)^{1+o(1)}, n\hat{\tau} - \alpha n^{2/3}]\} = O(n^{-K}),$$

for all $K > 0$. Here $T = T(u)$ is the running time for a fixed initial vertex $u$. That is, switching to $C_n(u)$ and using the union bound,

$$P\left\{\bigcap_{u \in [n]} \{C_n(u) \sim n[1-(p\pi+q)^d] \text{ or } C_n(u) = O(\omega^{-2} n^{2/3}(\ln n)^{1+o(1)})\}\right\}$$
(6.20)
$$\geq 1 - O(n^{-K+1}).$$

*... and existence of a giant component.* Corollary 1 implies that whp $G_{np^*}$ contains a component of size $\omega^{-1} n^{2/3}$ at least. As the multigraphs $G_{np^*}$ and $G_{np}$, $p > p^*$, can easily be coupled in such a way that $G_{np^*} \subset G_{np}$, we conclude that whp $G_{np}$ also has such a component. Since

$$\omega^{-1} n^{2/3} \gg \omega^{-3/2} n^{2/3},$$

it follows from (6.20) that whp this component has to be giant, of size close to $n[1-(p\pi+q)^d]$.

**Acknowledgments.** I owe a debt of gratitude to Itai Benjamini for posing the problem and for many stimulating discussions. The joint work with Jozsi Balogh on bootstrap site percolation gave me hope that viewing edge percolation as a dynamic process might help to alleviate enumerational difficulties which appear unavoidable in a static approach. I thank a referee for helpful comments.

# REFERENCES

[1] ALDOUS, D. J. and PITTEL, B. (2000). On a random graph with immigrating vertices: Emergence of the giant component. *Random Structures Algorithms* **17** 79–102. MR1774745

[2] ALON, N., BENJAMINI, I. and STACEY, A. (2004). Percolation on finite graphs and isoperimetric inequalities. *Ann. Probab.* **32** 1727–1745. MR2073175

[3] ARONSON, J., FRIEZE, A. and PITTEL, B. G. (1998). Maximum matchings in sparse random graphs: Karp–Sipser revisited. *Random Structures Algorithms* **12** 111–177. MR1637403

[4] BALOGH, J. and PITTEL, B. (2007). Bootstrap percolation on a random regular graph. *Random Structures Algorithms* **30** 257–286. MR2283230

[5] BENDER, E. A. and CANFIELD, E. R. (1978). The asymptotic number of labelled graphs with given degree sequences. *J. Combin. Theory Ser. A* **24** 296–307. MR0505796

[6] BENJAMINI, I. (2002). Private communication.

[7] BOLLOBÁS, B. (1980). A probabilistic proof of an asymptotic formula for the number of labelled regular graphs. *European J. Combin.* **1** 311–316. MR0595929

[8] BOLLOBÁS, B. (1982). The asymptotic number of unlabelled regular graphs. *J. London Math. Soc.* **89** 201–206. MR0675164

[9] BOLLOBÁS, B. (1984). The evolution of random graphs. *Trans. Amer. Math. Soc.* **286** 257–274. MR0756039

[10] BOLLOBÁS, B. (2001). *Random Graphs*, 2nd ed. Cambridge Univ. Press. MR1864966

[11] BORGS, C., CHAYES, J. T., VAN DER HOFSTAD, R., SLADE, G. and SPENCER, J. (2005). Random subgraphs of finite graphs. I. The scaling window under the triangle condition. *Random Structures Algorithms* **37** 137–184. MR2155704

[12] BORGS, C., CHAYES, J. T., VAN DER HOFSTAD, R., SLADE, G. and SPENCER, J. (2005). Random subgraphs of finite graphs. II. The lace expansion and the triangle condition. *Ann. Probab.* **33** 1886–1944. MR2165583

[13] DURRETT, R. (1985). Some general results concerning the critical exponents of percolation processes. *Z. Wahrsch. Verw. Gebiete* **69** 421–437. MR0787606

[14] DURRETT, R. (2005). *Probability*: *Theory and Examples*, 3rd ed. Wadsworth and Brooks/Cole, Pacific Grove, CA. MR1068527

[15] DURRETT, R. (2006). *Random Graph Dynamics*. Cambridge Univ. Press. MR2271734

[16] ERDŐS, P. and RÉNYI, A. (1960). On the evolution of random graphs. *Publ. Math. Inst. Hungar. Acad. Sci.* **5** 17–61. MR0125031

[17] GRIMMETT, G. (1989). *Percolation*. Springer, New York. MR0995460

[18] JANSON, S., ŁUCZAK, T., KNUTH, D. E. and PITTEL, B. (1993). The birth of the giant component. *Random Structures Algorithms* **4** 233–358. MR1220220

[19] ŁUCZAK, T. (1990). Component behavior near the critical point of the random graph process. *Random Structures Algorithms* **1** 287–310. MR1099794

[20] ŁUCZAK, T. (1990). On the number of sparse connected graphs. *Random Structures Algorithms* **1** 171–174. MR1138422

[21] ŁUCZAK, T., PITTEL, B. and WIERMAN, J. C. (1994). The structure of random graph near the point of the phase transition. *Trans. Amer. Math. Soc.* **341** 721–748. MR1138950

[22] MOLLOY, M. and REED, B. (1995). A critical point for random graphs with a given degree sequence. *Random Structures Algorithms* **6** 161–179. MR1370952

[23] MOLLOY, M. and REED, B. (1998). The size of the giant component of a random graph with a given degree sequence. *Combin. Probab. Comput.* **7** 295–305. MR1664335

[24] MOON, J. W. (1970). *Counting Labelled Trees*. Canadian Mathematics Congress, Montreal. MR0274333
[25] PITTEL, B. (2001). On the largest component of the random graph at a near critical stage. *J. Combin. Theory Ser. B* **82** 237–269. MR1842114
[26] PITTEL, B., SPENCER, J. and WORMALD, N. (1996). Sudden emergence of a giant $k$-core in a random graph. *J. Combin. Theory Ser. B* **67** 111–151. MR1385386
[27] WORMALD, N. C. (1995). Differential equations for random processes and random graphs. *Ann. Appl. Probab.* **5** 1217–1235. MR1384372

DEPARTMENT OF MATHEMATICS
OHIO STATE UNIVERSITY
231 W. 18TH AVENUE
COLUMBUS, OHIO 43210
USA
E-MAIL: bgp@math.ohio-state.edu